\documentclass[a4paper,12pt]{article}
\usepackage[cp1251]{inputenc}
\usepackage{graphicx}
\usepackage{amssymb,amsfonts}
\newtheorem{remark}{Remark}%
\newtheorem{theorem}{Theorem}
\newtheorem{proposition}{Proposition}
\newtheorem{lemma}{Lemma}
\newtheorem{definition}{Definition}

\def\bec{\begin{cor}}
\def\bedef{\begin{dfn}}
\def\bel{\begin{lem}}
\def\beprop{\begin{prop}}
\def\bet{\begin{thm}}

\def\enc{\end{cor}}
\def\endef{\end{dfn}}
\def\endprf{\noindent$\Box$}
\def\enl{\end{lem}}
\def\enprop{\end{prop}}
\def\ent{\end{thm}}

\def\prl{\noindent{\it Lemma proof. }}
\def\prt{\noindent{\bf Theorem proof. }}
\def\prp{\noindent{\bf Proposition proof. }}

\newcommand{\NN}{\text{\rm \hbox{I\kern-.2em\hbox{N}}}}

\title{Order positive  fields $\rm II$}
\date{}
\author{Margarita Korovina
 and Oleg Kudinov
\footnote{
M.V. Korovina's research was  carried out in the framework of the State Task to the
A.P. Ershov Institute of Informatics Systems (Project no. FWNU-2021-0003).
O.V. Kudinov's research was carried out in the framework of the State Task to the Sobolev Institute of Mathematics (Project FWNF–2022–0011)
and supported by project RNF 23-11-00170.
}
}
\begin{document}
\maketitle
\begin{abstract}
 This paper is a part of ongoing research on order positive fields started in \cite{korkud_al_1}.
%
We prove that the real closure of an order positive field even in non-Archimedean case is also order positive.

\end{abstract}

\section{Introduction}
In this paper  we continue to investigate the class of order positive fields introduced in \cite{korkud_al_1}.
Let us recall that this class is a proper extension of the class of computable ordered fields and differs from them by  the absence of the requirement of decidability of equality.

  The field of the primitive recursive real numbers and extensions of the rational numbers by computable sequences of computable real numbers are some of the examples of such fields. Note that in such cases equality is not decidable in general.
A significant role in the theory of computable fields is played by the well-known Ershov-Madison theorem \cite{Ersh_PRCM,Madison} which provides tools for effective construction of a computable presentation of the real closure of a given computable field.
Therefore the question whether the generalised result holds in order positive fields naturally arises.
In this paper we give a positive answer to this question.
For this, we give a characterization of the theory of finite extensions of order positive fields and on this basis establish that the class of order positive fields is closed under the real closure construction.

\section{Basic notions and definitions }

We refer the reader to \cite{vandervarden,Leng} for the notions and properties of polynomials and real closed fields,
to \cite{Viggo}
 for  basic definitions and fundamental concepts of computable rings and field,
 to \cite{Ersh_PRCM,Goncharov} for computable model theory, to \cite{Ersh_TN} for numbering theory and to \cite{korkud_al_1} for properties and examples of order positive fields.
 We fix a standard computable numbering $q:\omega\to \mathbb{Q}$  of the rational numbers and when it is clear from a context we use  the notation $q_n$ for $q(n)$,
 the standard notation $c_n:\omega^n\to\omega$ for Cantor's numbering of $n$--tuples. We also use the notations
  $[A,B]_L=\{x\in L\mid A\leq x\leq B\}$ and $(A,B)_L=\{x\in L\mid A<x < B\}$ for closed and open intervals in an ordered field $L$.

 In this paper we consider  a countable  ordered field $(F,\nu)$
with its numbering $\nu:\omega\to F$ as a $\sigma$-structure in the language $\sigma=(+,\cdot,0,1,<)$.
Let us recall the notion of an order positive field and their properties from \cite{korkud_al_1}, which we will use in this paper.
\begin{definition}\label{d_opf}(\cite{korkud_al_1})
A field $(F,\nu)$ is called order positive and the numbering $\nu$ is called its order positive presentation  if the following conditions hold:
\begin{enumerate}
 \item[1.] $(F,\nu)$ is an effective algebra, i.e., all operations are computable.

\item[2.] $\nu^{-1}(<)=\{(n,m)\mid \nu(n)<\nu(m)\}$  is computably enumerable.

  \item[3.] There exists a partial computable function $g:\omega\to\omega$ such that
\[
\nu(n)\neq 0\rightarrow g(n)\downarrow \wedge \, \nu(n)^{-1}=\nu(g(n)).
\]
\end{enumerate}
\end{definition}

\begin{proposition}\label{prop_1}(\cite{korkud_al_1})
If $(F,\alpha)$ is an effective algebra and $\alpha^{-1}(<)$ is computably enumerable then one can effectively construct a numbering $\beta$ induced by $\alpha$ such that
$(F,\beta)$ satisfies the conditions (1--3).
\end{proposition}
Below, in all proofs of order positivity we are going to use Proposition~\ref{prop_1} that allows us to avoid the verification  of the existence of  a function $g$
from Definition~\ref{d_opf}.
Examples of order positive fields including the primitive recursive real numbers can be found in \cite{korkud_al_1}.

By analogy with \cite{frohlich}
 we introduce the notion of an index of an order positive field that allows us to consider computable sequences
 $\{F_i\}_{i\in\omega}$ of order positive fields in a standard way, by requiring computability of sequences of certain indices.
\begin{definition}
Let $(F,\nu)$  be an order positive field. The number $i_F=c_4(a,b,c,d)$ is called its index if the following equalities hold
for all $k_1,\, k_2, \, n\in \omega$:
\begin{itemize}
\item $\varkappa_{a}=g_1$ and  $\nu(k_1)+\nu(k_2)=\nu(g_1(k_1,k_2))$,
\item $\varkappa_{b}=g_2$ and $\nu(k_1)\cdot\nu(k_2)=\nu(g_2(k_1,k_2))$,
\item $\varphi_{c}=g_3$ and $\nu(n)^{-1}=\nu(g_3(n))$ if  $\nu(n)\neq 0$,
\item $W_d=\{c_2(k_1,k_2)\mid \nu(k_1)<\nu(k_2)\}$,
\end{itemize}
where $\varkappa$ is  Kleene's numbering of the 2-arity patrial computable functions,  $\varphi$ is Kleene's numbering of the 1-arity partial computable functions,
$\{W_n\mid n\in\omega\}$ is Kleene's numbering of the computably enumerable sets.

\end{definition}
\section{Topological consistency}
Let us first say a few words about order topology, effectively open sets in this topology and their indices.
A basis of
the order topology on an ordered structure
$\mathcal{M}$ is the set of all open intervals
$(a,b)$,
where $a<b$.
If this basis is countable and has a numbering
$\lambda$
then a set $\cup_{n\in W_k} \lambda(n)$ is called  effectively open  with respect to this numbering for an appropriate computably enumerable set
 $W_k$ and an index of this effectively open set is $k$.

 It is clear that for an order positive field
$(K,\mu)$,
a basis of the order topology is the set of all intervals
$(\mu(n),\mu(m))_K$, where
$\mu(n)<\mu(m)$, and it has a natural numbering
$\lambda$.

\begin{proposition}\label{top}
Let $F$ be an ordered field
and $K$ be its algebraic extension.
Then,  $F$ is a subspace of $K$  in the order topology.
If $(F,\nu)$ and   $(K,\mu)$ are order positive fields, additionally
 $\nu\leq_c \mu$, and the set $\mathcal{U}\subseteq K$ is effectively open in the topology $\tau_K$ then
 the set $\mathcal{U}\cap F$ is effectively open in the topology  $\tau_F$,  i.e.,  effective versions  of the topologies are also consistent.
\end{proposition}
\prp
Let
$(a,b)_K$
 be an open interval in  $K$
and $(a,b)_K\cap F\neq \emptyset$.
We choose  $c\in F$ such that $c\in (a,b)_K$  then $0\in (a-c,b-c)_K$.
 Let $\alpha=b-c\in K$
then $\alpha>0$ and from algebraicity  (see \cite{Leng}) over  $F$ we have a bound $|\alpha^{-1}|<d$, where $d\in F$,
 that implies   $0<d^{-1}<\alpha$.
 By analogy for  $\beta=a-c<0$  we have  $|\beta^{-1}|<e$, where $e\in F$,  that
implies   $0<e^{-1}<-\beta$.
We get
  $(-e^{-1},d^{-1})_F\subset (\beta,\alpha)_{K}$ and all that remains is to shift the interval by  $c$
 then $c\in (c-e^{-1}, c+ d^{-1})_F\subset (a,b)_{K}$.
It follows from here that
$F$ is a subspace of  $K$, i.e.,
  an open  in $K$ set intersected with $F$ is an  open set in $F$.

Now we  show that if
   $\mathcal{U}=(a,b)_K=(\mu(n),\mu(m))_K$ then
$\mathcal{U}\cap F$ is effectively open in $F$.
Let $\nu(k)=\mu(h(k))$  for all $k$, where $h$ is a computable function.
Without having a uniform way to determine minimal polynomial for
 $\alpha=b-c$ we simply use an enumeration as follows: first $c=\nu(k_0)$
  for $k_0$ with the condition $a<\mu(h(k_0))<b$, then $d=\nu(k_1)=\mu(h(k_1))$
  for the number $k_1$  with the condition $0<\alpha^{-1}<d$ and $e=\nu(k_2)$
 with the condition $0<e^{-1}<\beta=c-a$.
In the computation above we have the computably enumerable set  $I=\{k_0\mid \mu(n)<\mu(h(k_0))<\mu(m)\}$
  and some partial computable functions $g_1(n,m,k_0)$ and $g_2(n,m,k_0)$
  with the domains $\{(n,m)\mid \mu(n)<\mu(m)\}\times I$
 that by $k_0$ compute indices $k_1$
   and $k_2$.

Finally, the set
   \[ (a,b)_K\cap F=\bigcup_{k_0\in I}\Big(\nu(k_0)-\nu(g_2(n,m,k_0)^{-1}),\nu(k_0)+\nu(g_1(n,m,k_0)^{-1} )\Big)_{F}
   \]

\noindent is effectively open in $F$ according to  order positivity of $\nu$.

From this,  the transition from
  $\mathcal{U}$
 to $\mathcal{U}\cap F$  follows in the general case.
\endprf
\section{ The real closure of an order positive field}
Let    $\overline{F}=[F]_{rcl}$  denote  the real closure of $F$,  $D_f $ denote the discriminant of a polynomial  $f$,
  ${\rm Res}(f,g)$ denote the resultant of polynomials  $f$ and $g$.

Below we often refer to a construction of intermediate order positive presentations of fields  $F\leq F_1\leq \overline{F}$,
where the extension $F_1\geq F$ is finite. For example, an extension can be obtained by adding to $F$ one or few roots in $\overline{F}$ of a polynomial  $f\in F[x]$.
In such cases, we choose a convenient copy $\widetilde{F}$ of the field $F_1$ such that there exists  an isomorphism $\varphi:\widetilde{F}\to F_1$
and propose its order positive presentation $\widetilde{\nu}$.
 When needed, an  order positive presentation can be constructed for the field $F_1$ as the composition $\varphi\circ \widetilde{\nu}$
 however
 this is not always necessary.

Moreover, it is worth noting that the number of all roots of a polynomial  $f\in F[x]$ in $\overline{F}$ and likewise the number of all its roots in  $[A,B]_{\overline{F}}$, do not depend on the choice of $\overline{F}$, as these numbers are uniquely defined by the Sturm's method  \cite{vandervarden}  although the latter requires the equality verification for elements of $F$, i.e., it is  not quite constructive.

 For numbered field and their extensions  we use the following notation. For $\widetilde{F}\geq F$ we write $(\widetilde{F},\widetilde{\nu})\succeq (F,\nu)$ if $\nu$ is reducible to  $\widetilde{\nu}$,
i.e., for  an appropriate computable function $\nu= \widetilde{\nu}\circ f$.
\begin{proposition}\label{prop}
Let $(F,\nu)$  be an order positive field.
For every $n>0$, a unitary polynomial $f\in F[x]$ with   ${\rm deg}(f)=n$ and for all $A,\, B\in F$
with the condition $A<B$ the following statements  hold:
\begin{enumerate}
\item The sets \[\{y\in F\mid y>\sup_{x\in[A,B]_{\overline{F}}}f(x)\}, \] \[\{y\in F\mid y<\sup_{x\in[A,B]_{\overline{F}}}f(x)\}\]
are effectively open and their indices can be effectively found  by indices of $f,\,A,\,B$ and  $F$.

\item If  $f$ does not have multiple roots then for an appropriate finite order positive extension  $(\widetilde{F},\widetilde{\nu})\succeq (F,\nu)$
    one can effectively construct the set
   $\{\langle k_i,m_i\rangle\in\omega^2\mid 0<i\leq m\}$
   of
 indices of intervals  $(A_i,B_i)_{\overline{F}}$,  where $A_i=\widetilde{\nu}(k_i)$ and $B_i=\widetilde{\nu}(m_i)$, such that
  for any $i<m$
  either $A_{i+1}>B_i$
 or $k_{i+1}=m_i$
 and all roots $f$ from $\overline{F}$
 belong to $\bigcup_{i\leq m}(A_i,B_i)_{\overline{F}}$,
 in every $(A_i,B_i)_{\overline{F}}$
  there is exactly one root of   $f$.
An index of $\widetilde{F}$ and an index of a function reducing $\nu$ to $\widetilde{\nu}$ can be effectively found using indices of $f$ and $F$.

  \item If  $f$  does not have multiple roots then the field $\widetilde{F}$ obtained by joining all roots of $f$
  from $\overline{F}$ to $F$ is order positive, it has an order positive presentation $\widetilde{\nu}$ such that   $\nu<\widetilde{\nu}$.
An index  of $\widetilde{F}$, an index of a function reducing $\nu$ to $\widetilde{\nu}$
 and a finite tuple containing  exactly one index of each root of
 $f$ in $\widetilde{F}$
      can be effectively found by  indices of $f$ and $F$.

\end{enumerate}

\end{proposition}
\prp
We simultaneously prove all statements by  induction with the parameter $n={\rm deg}(f)$.
For that, we use the following definition and lemmas.

For a polynomial $f$ let us define $f_\delta(x)=f(x)-\delta\cdot x$, where $\delta\in F$.
We say that an element $\delta\in F$ is {\em admissible} if $f_\delta^\prime(x)=f^\prime(x)-\delta$ does not have multiple roots.
From  the definition  above it follows that:
\begin{itemize}
\item[a)] almost all elements of the field $F$, (apart from a finite number)  are admissible,
\item[b)] the set
 $\{n\in\omega\mid \delta=\nu(n) \mbox{ is  admissible}\}$ is computably enumerable.
\end{itemize}
The item a)  follows from the definition, i.e.,
an element is non-admissible if it is a root of the  discriminant $D_{f^\prime(x)-\delta}$.
For verifying the item b) it is sufficient to enumerate
$n\in\omega$ such that the discriminant of the polynomial
$f^\prime(x)-\delta$
differs from 0. This is possible by assumption that $(F,\nu)$ is order positive.
Moreover the item b) holds uniformly on $f$, i.e.,
after the construction of a numbering of the polynomial ring $F[x]$ (see\cite{frohlich})  one can compute an index of this set
using an index of $f$.
\begin{lemma}\label{lemma_1}
For every   $f\in F[x]$,   $A,\, B\in F$
with $A<B$  and $y\in F$
the following statements  hold
\begin{enumerate}
\item
 $y>\sup_{x\in[A,B]_{\overline{F}}}f(x)$ if and only if there exists admissible
 $\delta>0$ such that
 $y>\sup_{x\in[A,B]_{\overline{F}}}f_\delta(x)+ \delta\cdot\max\{|A|,\, |B|\}$.
 \item
 $y<\sup_{x\in[A,B]_{\overline{F}}}f(x)$ if and only if there exists admissible
 $\delta>0$ such that
 $y<\sup_{x\in[A,B]_{\overline{F}}}f_\delta(x)- \delta\cdot\max\{|A|,\, |B|\}$.
 \end{enumerate}
\end{lemma}
\prl
It is worth noting that there exists
 $z=\sup_{x\in[A,B]_{\overline{F}}}f(x)=\max_{x\in[A,B]_{\overline{F}}}f(x)$ in $\overline{F}$
(the corresponding fact in $\mathbb{R}$ is an elementary property therefore it holds on all real closed fields).
This remark is standard and we use similar transfer principle in  further  proofs.

\noindent $1.\Rightarrow)$
Since
$y\in(z,+\infty)_{\overline{F}}$
by the definition \ref{top}
we have
$(y-3\epsilon,y+3\epsilon)_F\subseteq (z,+\infty)_{\overline{F}} $ for appropriate  $\epsilon>0,\, \epsilon\in F$.
Then $y>z+2\epsilon$.
We choose $\delta>0$ in $F$ such that $\delta\cdot\max\{|A|,\,|B|\}<\epsilon$
and
$\delta$  is admissible.
Such $\delta$ is the required number since
 $y>\sup_{x\in[A,B]_{\overline{F}}}f_\delta(x)+ \epsilon$
by $|f(x)-f_\delta(x)|\leq\delta\cdot\max\{|A|,\,|B|\}$
for $x\in[A,B]_{\overline{F}}$ and  $y>\sup_{x\in[A,B]_{\overline{F}}}f(x)+2\epsilon$ and also
$\delta\cdot\max\{|A|,\,|B|\}<\epsilon$.

\noindent $1.\Leftarrow)$, $2.\Leftarrow)$
 By the inequality $|f(x)-f_\delta(x)|\leq \delta\cdot\max\{|A|,\,|B|\}$  for   ${x\in[A,B]_{\overline{F}}}$
the statement is straightforward.
\endprf
\begin{lemma} Let
$h_1(x),\dots,h_s(x)$
 be a finite collection of  $F[x]\setminus\{0\}$.
Then the set
\begin{eqnarray*}
&C=\{\alpha\in F\mid \mbox{ for } a(x)=f(x)-\alpha,\mbox{  and all }i=1,\dots,s \\
& {\rm Res}(a(x),h_i(x))\neq 0 \}
 \end{eqnarray*}
 is co-finite, in particular it is dense in  $F$.
\end{lemma}
\prl
For every  $i=1,\dots,s$ the equality ${\rm Res}(f(x)-\alpha,h_i(x))= 0$ holds only for a finite number of
 values $\alpha$, this implies the statement.
\endprf

\noindent {\bf Statement 1.}
We search for parts of required sets
 $X=\{y\in F\mid y>\sup_{x\in[A,B]_{\overline{F}}}f(x)\}$ and
$Y=\{y\in F\mid y<\sup_{x\in[A,B]_{\overline{F}}}f(x)\}$
using various finite extensions of the field
 $F$
such that the sets  $X$ and $Y$
will be the unions of the corresponding parts.
 An advantage of admissible
 $\delta\in F$ is that
 using inductive assumption by Statement 3 we can add all roots of the polynomial $f_\delta^\prime$ to $F$ and
 get an order positive extension   $\widetilde{F}$ by computing its index.
Below we describe our construction  of the sets $X$ and $Y$ in  details.
For every admissible $\delta>0,\, \delta\in F$, let us define the finite set $J_\delta$ (that could be empty) consisting of
all roots of $f_\delta^\prime$ and an appropriate order positive extension $\widetilde{F}_\delta$ containing  $J_\delta$ such that
its index is effectively depends on $m\in \nu^{-1}(\delta)$. We assume that from now $m$  and $\widetilde{F}_\delta$ are fixed.

Further, for every admissible $\delta>0, \delta\in F$,  we define the following sets:
\begin{eqnarray*}
& X_\delta=&\{y\in F\mid y>f_\delta(A)+ \delta\cdot\max\{|A|,\, |B|\}\mbox{ and }\\
&&y>f_\delta(B)+ \delta\cdot\max\{|A|,\, |B|\} \mbox{ and} \\
&& (\forall a\in J_\delta)\,  y>f_\delta(a)+ \delta\cdot\max\{|A|,\, |B|\} \mbox{in } \widetilde{F}_\delta\}
\end{eqnarray*}
and
\begin{eqnarray*}
& Y_\delta=&\{y\in F\mid y<f_\delta(A)- \delta\cdot\max\{|A|,\, |B|\}\mbox{ or }\\
&&y<f_\delta(B)- \delta\cdot\max\{|A|,\, |B|\} \mbox{ or} \\
&& (\exists a\in J_\delta)\,  y<f_\delta(a)- \delta\cdot\max\{|A|,\, |B|\} \mbox{ in } \widetilde{F}_\delta\}.
\end{eqnarray*}
It is easy to see that
 $X_\delta$ and $Y_\delta$ are intersections of effectively open in $ \widetilde{F}_\delta$ sets with $F$, so
 by Proposition~\ref{top} the sets $X_\delta$ and $Y_\delta$ are effectively open in $F$.
 Since the polynomial  $f_\delta$ attains its maximums in  $[A,B]_{\overline{F}}$ we have the following equivalences:
 \[ y\in X_\delta \leftrightarrow y>\sup_{x\in[A,B]_{\overline{F}}}f_\delta(x)+ \delta\cdot\max\{|A|,\, |B|\}\]
and
  \[ y\in Y_\delta \leftrightarrow y<\sup_{x\in[A,B]_{\overline{F}}}f_\delta(x)- \delta\cdot\max\{|A|,\, |B|\}.\]
By Lemma \ref{lemma_1} it follows that
\[ y\in X \leftrightarrow  \mbox{  there exists admissible }\delta>0\,\mbox{ such that } y\in X_\delta,\]
\[ y\in Y \leftrightarrow  \mbox{  there exists admissible }\delta>0\,\mbox{ such that } y\in Y_\delta,\]
so by Lemma 2 the sets  $X$  and $Y$ are effectively open in $F$.

\noindent  {\bf Statement 2.}
 A preliminary interval $[A,B]_{\overline{F}}$, $A,\, B \in \overline{F}$, that contains all roots of a polynomial $f=\sum_{i=0}^n a_ix^i$
from $\overline{F}$ can be found due to the well known estimation of its roots $\beta$:
\[\beta\leq \max\Big(1, \sum_{i=0}^n |a_i| \Big), \mbox{ at the same time } f(A)\neq 0,\, f(B)\neq 0.\]

 Since $f^\prime$ could have multiple roots we can not
straightforwardly use the inductive assumption
to get $\widetilde{F}$  joining the roots of $f^\prime$ from $\overline{F}$.
Instead using Lemma~2 we search for an admissible $\pm \delta\in F$ with $0<\delta<\lambda$ (the choice of $\lambda$ will be explained later) and polynomials $\widetilde{f}=f^\prime\pm\delta$
with the following requirements:
 \begin{itemize}
 \item  $\widetilde{f}$  does not have multiple roots,
 \item $|f^\prime(x)-\widetilde{f}(x)|<\frac{\lambda}{2}$ for every $x\in[A,B]_{\overline{F}}$,
  \item $|f^\prime(\beta)|>\lambda$ for every $\beta\in[A,B]_{\overline{F}}$ with $f(\beta)=0$,
 \item $f^\prime(A),\, f^\prime(B)\not\in \{-\delta,+\delta\}$.
 \end{itemize}
It is worth noting that $\lambda$ can be effectively found using the decomposition in $F[x]$
(see \cite{vandervarden}):

 \[
  f(x)C(x)+f^\prime(x) D(x)={\rm Res}(f,f^\prime) \mbox{ for } C(x),\, D(x)\in F[x],
 \]
where
 ${\rm deg}(D)< {\rm deg}(f)=n$,
and the bounds:
 \[
 |f^\prime(\beta)|>\frac{|{\rm Res}(f,f^\prime)|}{\sup_{x\in[A,B]_{\overline{F}}}|D(x)|}>\lambda.
 \]
If $D(x)=\sum^m_{i=0}\beta_ix^i$, $m<n$,
(while we do not know  whether coefficients are equal to zero or not)
then the following bound
\[
|D(x)|<\sum^m_{i=0}(|\beta_i|+1)E^i+1, \mbox{ where } E= max(|A|,\, |B|),
\]
is straightforward
for any
$x\in[A,B]_{\overline{F}}$.
This means that we can bound  $D$ from above
 on  $[A,B]_{\overline{F}}$, which leads to a successful search of a suitable   $\lambda \in F$.

Let the finite set  $J_\delta$ be the union of $A,\,B$ and all roots of the polynomial $\widetilde{f}$ for a chosen $\delta$.
It is obvious that its elements are not roots of the polynomial $f$.
The construction below is intended to put the roots of the polynomial $f$ and certain pairs of adjacent elements of $J_\delta$ into
a bijective correspondence.

Assume that $\beta$ is a root of $f$ on  $[A,B]_{\overline{F}}$.
Then for some pairs
$\langle a,\,b\rangle\in {\overline{F}}\times {\overline{F}}$ the following holds
 \begin{itemize}
 \item $a<b$,  $a,\, b\in [A,B]_{\overline{F}}$  and $\beta\in (a,b)_{\overline{F}}$,
 \item $f^\prime(a)=0$ or $a=A$,
 \item $f^\prime(b)=0$  or $b=B$,
 \item $f(a)\cdot f(b)<0$.
 \end{itemize}
 From these pairs we choose a pair $\langle a, b\rangle$  of adjacent elements such that $a$ is the maximum  and $b$ is the minimum  of all possible elements.

 \noindent Let us consider all  cases for $a$.
  Assume $f^\prime(a)=0$ and $a \geq A$.
 Then by  the intermediate value theorem there exists
$\beta_1\in(a,\beta)_{\overline{F}}$
such that
 $f^\prime(\beta_1)\in \{-\delta,\delta\}$,
moreover
$f^\prime(\beta_1)\cdot f^\prime(\beta)>0$.

Let us choose $\beta_1$ that is  the maximum of them intersected with  $ J_\delta$.
At the same time
$f$ and $f^\prime$ do not change signs in the interval $(\beta_1,\beta)_{\overline{F}}$.
So $\beta_1$ is the left element of a required pair.
\noindent Let $a = A$  and $f^\prime(a)\neq 0$.
In this case, an element $z\in \overline{F}$ may also be found such that
$z\in (a,\beta)$ and
$f^\prime(z)=\pm \delta$
(the sign of $f^\prime(z)$  is the same as one of
$f^\prime(\beta)$).
If such elements exist then we put  $\beta_1=z$, where $z$ is the maximum of them intersected with  $ J_\delta$, otherwise $\beta_1=A$.
In both cases $f$ does not change the  sign in the interval
 $(\beta_1,\beta)_{\overline{F}}$. So $\beta_1$ is the left element of a required pair.

 \noindent
We consider similarly the cases for  $b$ and choose $\beta_2\in J_\delta$ that is the right element of a required pair.

We already established the following: if $f(\beta)=0$ and $\beta\in(A,B)_{\overline{F}}$  then for
 the adjacent elements $\beta_1<\beta_2$ of $J_\delta$,
we have
$\beta\in(\beta_1,\beta_2)_{\overline{F}}$
and
$f(\beta_1)\cdot f(\beta_2)<0$.
On the other hand if there exist adjacent elements $\beta_1,\,\beta_2\in J_\delta$ with the condition $f(\beta_1)\cdot f(\beta_2)<0$
then by the intermediate value theorem there exists an   element
 $\beta\in\overline{F}$
such that $\beta\in(\beta_1,\beta_2)_{\overline{F}}$ and $f(\beta)=0$.
There cannot be more than one such element. Indeed,
assume contrary that there exists another one
 $\beta^\ast\in (\beta_1,\beta_2)_{\overline{F}}$ different from $\beta$
such that $f(\beta^\ast)=0$.
Then by Rolle's theorem \cite{vandervarden} there exists $a^\ast\in (\beta^\ast,\beta)_{\overline{F}}$ such that
$f^\prime(a^\ast)=0$. This contradicts  the choice of the elements $a,\,\beta_1$, namely $a\leq \beta_1<\beta$.
Let us recall that
$a$ is the maximum, at the same time each root of  $f$ is located between adjacent elements $\beta_1<\beta_2$ of $J_\delta$ such that $f(\beta_1)f(\beta_2)<0$. Then
 $ \widetilde{F}=F(\{\alpha\mid\alpha\in J_\delta\})$.
Above we already noted that for
$\alpha\in J_\delta$  we have $f(\alpha)\neq 0$.
Therefore we choose previously constructed intervals
$(\beta_1,\beta_2)_{\overline{F}}$
 (with the condition $f(\beta_1)f(\beta_2)<0$)
 as required ones
 $(A_i,B_i)_{\overline{F}}$.
A reduction function from  $\nu$ to $\widetilde{\nu}$ will be constructed below.

\noindent  {\bf Statement 3.}
Let
$(A_i,B_i)_{\overline{F}}$, $i=1,\dots,m$ be intervals constructed in Statement 2.
If $\alpha_i$ is a roof of $f$ from $(A_i,B_i)_{\overline{F}}$ then it is sufficient to represent $\widetilde{F}$ as
$F(\alpha_1,\dots,\alpha_{m-1})(\alpha_m)$ and  explain how to construct an order positive presentation of the field
$F_i=F_{i-1}(\alpha_{i})$, i.e., show how to add  one root at a time to an already constructed field.
Formally,
$F_i=F(\alpha_1,\dots,\alpha_i)$,
$F_{i+1}=F_i(\alpha_{i+1})$,
$i< m$,
and thus
$F_m=\widetilde{F}$ will be constructed.

Without loss of generality we assume that
 $\alpha$ is the unique  root of $f$ in
the interval
$(A,B)_{\overline{F}}$.
We are going to construct the diagram of the field $F(\alpha)$ in the language $\sigma=(+,\cdot,0,1,<)$.
It is clear that for that it is sufficient to computably enumerate all indices of polynomials from
$\{a(x)\in F[x]\mid a(\alpha)>0\}$
since every element of $F(\alpha)$ can be represented as $a(\alpha)$, where $a\in F[x]$ and  ${\rm deg}(a)< {\rm deg}(f)$.
Therefore the numbering $\mu$ of  the polynomials  $a\in F[x]$ leads to the numbering $\widetilde{\nu}$
of all elements of the field $F(\alpha)$
according $\widetilde{\nu}(n)=\mu(n)(\alpha)$ that induces an order positive presentation of $F(\alpha)$. In particular,
some index of  $\alpha$ is  a $\mu$-index of the polynomial $a(x)\equiv x$.
First in the ring $F(x)$ we divide every  polynomial $a\in F[x]$ by $f$ with a reminder $d$: $a(x)=f(x)\cdot b(x)+d(x)$,
$d(x)=\sum_{i=0}^{n-1}\alpha_ix^{n-1-i}$.
It is worth noting that while the elements
$\alpha_i$, $i\leq n-1$, can be found in a standard way by the division algorithm in $F[x]$ the exact degree of $d(x)$ is not known.
Suppose we found
 $\delta >0$ and $\epsilon>0$
in $F$, $k\leq n-1$,
such that

\[
\epsilon\cdot\sum_{i=0}^k c^{n-1-i}<\delta\wedge\bigwedge_{i=0}^k|\alpha_i|<\epsilon
\wedge \big(|\alpha_{k+1}|>\epsilon\vee k=n-1\big),
\]
where
$c=\max(|A|,|B|)$.
In this case we say that $(\delta,\epsilon,k)$ is a consistent 3-tuple and define
$\widetilde{d}=\sum_{i=k+1}^{n-1}\alpha_ix^{n-1-i}$.
For consistent 3-tuple  and $x\in[A,B]_{\overline{F}}$   $|d(x)-\widetilde{d}(x)|<\delta$ holds,
at the same time
either $\widetilde{d}\equiv 0$
or  $\alpha_{k+1}>\epsilon$.
Let us use the following notation.

\noindent {\bf Condition $\ast$ on $ (\delta,\epsilon,k)$}:
{\it For some
$y_0,\,y_1\in F$
such that $y_0>0,\,y_1>0$,
both polynomials
$\widetilde{d}(x)-y_0$ and $\widetilde{d}(x)-y_1$ do not have multiple roots,
moreover for appropriate
$x_0$,  $x_1$  with

 \begin{eqnarray*}
 \left \{
\begin{array}{r}
\widetilde{d}(x_0)=y_0\\
 \widetilde{d}(x_1)=y_1
\end{array}
\right.
\end{eqnarray*}
we have
$A<x_0<x_1<B$
(in a suitable order positive extension  $\widehat{F}>F$)
and at the same time
$\inf_{x\in[x_0,x_1]_{\overline{F}}}\widetilde{d}(x)>\delta$  and  $f(x_0)f(x_1)<0$.
}

If for some consistent 3-tuple $(\delta,\epsilon,k)$ the condition~$\ast$ holds then
$d(x)>0$ for $x\in[x_0,x_1]_{\overline{F}}$, in particular, for $\alpha\in [x_0,x_1]_{\overline{F}}$ we obtain
$d(\alpha)=a(\alpha)>0$.

Converse is also valid:
 $d(\alpha)=a(\alpha)>0$ implies  the existence of a consistent 3-tuple $(\delta,\epsilon,k)$ with the condition~$\ast$.
For a given $a(x)\in F[x]$  we can uniformly enumerate all consistent 3-tuples $(\delta,\epsilon,k)$  with the condition~$\ast$
so we can enumerate all those $a(x)\in F[x]$ for which there exist such 3-tuples, in particular $\{a(x)\mid a(\alpha)>0\}$.

Thereby we can computably enumerate $a(x)\in F[x]$ with the condition $a(\alpha)>0$ that leads to construction of an order positive presentation of $F(\alpha)$.
It is worth noting that from the proof it follows that there is possible to compute (with respect to $\widetilde{\nu}$)
some indices of all roots  $\alpha_1,\dots,\alpha_m$ of $f$  in $\widetilde{F}$.
\endprf

\begin{remark}
The following example shows that in  non-Archimedean  case of  Statement 2 of Proposition~\ref{prop} the construction of an appropriate extension $(\widetilde{F},\widetilde{\nu})$ of the field
$(F,\nu)$ is unavoidable.
Let $F=\mathbb{Q}(t)$, where $t>\mathbb{Q}^+$.
Then  in the field $F$ the roots of the polynomial $f(x)=(x^2-t)(x^3-t)$ cannot be separated, i.e.,
there is no an interval $(\alpha,\beta)$ with $\alpha,\beta\in F$ such that $\sqrt{t}\in(\alpha,\beta)$ however
$\sqrt[3]{t}\not\in(\alpha,\beta)$.
\end{remark}
\begin{theorem}
The real closure of an order positive field is also order positive.

\end{theorem}
\prt
We create the full list $\{f_i(x)\}_{i\in\omega}$ of unitary polynomials without multiple roots, i.e.,
with non-zero discriminant $D_{f_i}\neq 0$ and construct
 a field $F^\ast$ containing all roots from $\overline{F}$ of all polynomials $f_i$ and generated by them.
We construct the sequence $\{F_i\}_{i\in\omega}$ of order positive fields such that
 $F_{i+1}\supseteq F_i$ and
$F^\ast=\bigcup_{i\in\omega} F_i$.
This sequence is computable in the sense that indices of fields are computable.

Construction:
\begin{itemize}
\item $F_0=F$,
\item $F_{i+1}$ is obtained by adding to $F_i$  all roots of $f_i$ from $\overline{F}$,
$\nu_i\leq_{h_i}\nu_{i+1}$ (by Statement 3 in Proposition \ref{prop}),
\item moreover it is straightforward to ensure that $\{h_i\}_{i\in\omega}$ is a computable sequence.
\end{itemize}
By recursion, we compute the index of the function $H(i,j)$ that reduces  $\nu_i$ to $\nu_j$ for $i\leq j$.
Then
 $\nu^\ast(c(i,n))=\nu_i(n)$ is an order positive presentation of  $F^\ast$ that
  is isomorphic to $\overline{F}$ over $F$, i.e., the isomorphism is identical on $F$.
\endprf

Authors and Affiliations:

\medskip

M.V. Korovina,

A.P. Ershov Institute of Informatics Systems SB RAS,

6, Acad. Lavrentjev pr., Novosibirsk 630090, Russia

e-mail: rita.korovina@gmail.com

O.V. Kudinov,

Novosibirsk State University,

2, Pirogova str., Novosibirsk, 630090, Russia

 A.P. Sobolev Inst. Math. SB RAS

4, Acad. Koptuga pr.,  Novosibirsk 630090, Russia

 e-mail: kud@math.nsc.ru

\end{document}